\documentclass[11pt]{article}

\usepackage{tikz-cd}

\usepackage[utf8]{inputenc}
\usepackage[OT2]{fontenc}
\usepackage[T1]{fontenc}

\usepackage{palatino}

\usepackage[english]{babel}
\usepackage{comment}
\usepackage{amsmath}
\usepackage{amsfonts}
\usepackage{mathrsfs}
\usepackage{amssymb}
\usepackage{geometry}
\usepackage{stmaryrd}
\usepackage{chemfig}
\usepackage{mathtools}
\usepackage{amsthm}
\usepackage{hyperref}
\usepackage{wasysym}

\theoremstyle{plain}
\newtheorem{theo}{Theorem}[section]
\newtheorem{prop}[theo]{Proposition}
\newtheorem{coro}[theo]{Corollary}
\newtheorem{lemma}[theo]{Lemma}

\theoremstyle{definition}
\newtheorem{madef}[theo]{Definition}
\newtheorem{nota}[theo]{Notation}
\newtheorem*{nota*}{Notation}
\newtheorem{rmq}[theo]{Remark}
\newtheorem{ex}[theo]{Example}
\newtheorem{fact}[theo]{Fact}
\newtheorem{madefprop}[theo]{Definition/Proposition}

\theoremstyle{remark}

\newcommand{\un}{\hat{1}}
\newcommand{\zero}{\hat{0}}
\renewcommand{\L}{\mathcal{L}}
\newcommand{\G}{\mathcal{G}}

\newcommand{\Ind}{\mathrm{Ind}}
\newcommand{\A}{\mathcal{A}}

\newcommand{\rk}{\mathrm{rk}}
\renewcommand{\S}{\mathcal{S}}

\newcommand{\FY}{\mathbb{FY}}
\newcommand{\FYa}{\mathrm{FY}}
\newcommand{\I}{\mathcal{I}}
\newcommand{\aff}{\mathrm{aff}}

\newcommand{\wond}{\mathrm{wond}}

\newcommand{\At}{\mathrm{At}}

\newcommand{\totoreq}{\trianglelefteq}

\renewcommand{\d}{\mathrm{d}}

\newcommand{\N}{\mathbb{N}}
\newcommand{\Cx}{\mathbb{C}}
\newcommand{\Z}{\mathbb{Z}}

\newcommand{\ANM}{\mathsf{ANM}}
\newcommand{\ONM}{\mathsf{ONM}}

\newcommand{\B}{\mathcal{B}}

\newcommand{\EL}{\lambda_{\vartriangleleft}}

\renewcommand{\P}{\mathcal{P}}

\newcommand{\Circ}{\mathcal{C}}

\newcommand{\Supp}{\mathsf{Supp}}

\newcommand{\w}{\mathsf{w}}

\newcommand{\LBS}{\mathfrak{LBS}}

\newcommand{\Ne}{\mathcal{N}}
\newcommand{\Q}{\mathbb{Q}}
\newcommand{\weight}{\mathsf{w}}

\newcommand{\Tor}{\mathbb{T}}
\newcommand{\R}{\mathbb{R}}
\newcommand{\Mod}{\mathcal{M}}
\newcommand{\irr}{\mathrm{irr}}

\title{Supersolvability of built lattices and Koszulness of generalized Chow rings}
\author{Basile Coron}
\date{}

\begin{document}
\maketitle
\begin{abstract}
We give an explicit quadratic Gröbner basis for generalized Chow rings of supersolvable built lattices, with the help of the operadic structure on geometric lattices introduced in a previous article. This shows that the generalized Chow rings associated to minimal building sets of supersolvable lattices are Koszul. As another consequence, we get that the cohomology algebras of the components of the extended modular operad in genus $0$ are Koszul.\end{abstract}
\tableofcontents
\section{Introduction}
In \cite{FY_2004} Feichtner and Yuzvinsky defined algebras $\FYa(\L, \G)$ for every pair of a geometric lattice $\L$ and a building set $\G \subset \L$ (such datum $(\L, \G)$ will be called a built lattice). In the realizable case those algebras are the cohomology rings of the wonderful compactifications introduced by De Concini-Procesi \cite{de_concini_wonderful_1995}. If $\G$ is equal to $\L\setminus \{\zero\}$ we get the so-called combinatorial Chow ring of $\L$. Those rings are known to satisfy very strong properties such as Poincaré duality or even the Kähler package (see Adiprasito-Huh-Katz \cite{Huh_2018} for combinatorial Chow rings and Pagaria-Pezzoli \cite{Pagaria_2021} for general Feichtner--Yuzvinsky rings).\\

An important property of algebras which is still a largely open question for Feichtner--Yuzvinsky algebras is that of Koszulness. In plain English, Koszulness means that the algebra in question has a weight grading such that the algebra is generated by elements of weight 1, the relations between elements of weights 1 are generated by elements of weight 2 (i.e. the algebra is quadratic), the relations between relations are generated in weight 3 and so on. Koszulness is a particularly interesting property to ask of the cohomology ring of a formal space because it allows a direct computation of other rational homotopy invariants such as the rational homotopy Lie algebra (see Berglund \cite{Berglund_2014}). Since the wonderful compactifications of hyperplane arrangements are known to be formal, it is natural to ask which Feichtner--Yuzvinsky algebras are Koszul, a question raised by Dotsenko in \cite{Dotsenko_2022}. \\

A classical way to prove the Koszulness of a given algebra is to find a quadratic Gröbner basis for this algebra. Feichtner and Yuzvinsky computed explicit Gröbner bases for the Feichtner--Yuzvinsky rings, but those bases are almost never quadratic. In fact, the Feichtner--Yuzvinsky rings themselves are not necessarily quadratic. One of the first results proving the Koszulness of some Feichtner--Yuzvinsky algebras was given by Dotsenko who proved that the Feichtner--Yuzvinsky algebras associated to the complete graphs with building set of connected components are Koszul. To put it in a nutshell, Dotsenko introduced an explicit order on the generators of the Feichtner--Yuzvinsky rings and then used the operadic structure on this collection of rings to construct a bijection between the algebraic normal monomials associated to the latter order and relations of degree $2$, and the operadic normal monomials obtained in a previous work via Gröbner bases for operads (Dotsenko-Khoroshkin \cite{DK_2010}). By a dimension argument this implies that the relations of weight $2$ form a quadratic Gröbner basis of the Feichtner--Yuzvinsky rings in question. More recently, Mastroeni-McCullough \cite{MM_2022} proved that the combinatorial Chow rings are all Koszul, using the notion of Koszul filtrations. \\

In \cite{Stanley_1972} Stanley introduced the class of ``supersolvable'' lattices.
\begin{madef}[Stanley, \cite{Stanley_1972}]
A lattice $\L$ is called \textit{supersolvable} if it admits a maximal chain $\Delta$ such that for every chain $K$ in $\L$, the sublattice generated by $\Delta$ and $K$ is distributive. 
\end{madef}
Supersolvable lattices have very nice properties in general. In particular, we have the following classical result.
\begin{theo}[Yuzvinsky, \cite{Yuzvinsky_2001}]
The Orlik--Solomon algebra of a supersolvable lattice admits a quadratic Gröbner basis. 
\end{theo}

In this article we prove a similar result for Feichtner--Yuzvinsky algebras. We first introduce a notion of supersolvability for built lattices (which coïncides with the usual supersolvability when taking the maximal building set) and we prove the following theorem.
\begin{theo}\label{theomainintro}
Let $(\L, \G)$ be a supersolvable built lattice. The algebra $\FYa(\L, \G)$ admits a quadratic Gröbner basis. 
\end{theo}
In order to prove this result we will generalize the strategy of Dotsenko, by using the extended operadic structure introduced in $\cite{Coron_2022}$. \\

Theorem \ref{theomainintro} immediately implies that combinatorial Chow rings of supersolvable lattices have quadratic Gröbner bases, strengthening the result of Mastroeni-McCullough for supersolvable lattices. Turning our attention towards minimal building sets, Theorem \ref{theomainintro} also gives us the following result. 

\begin{theo}\label{theogminintro}
Let $\L$ be a supersolvable lattice and $\G_{\min}$ the building set of irreducible elements of~$\L$. The algebra $\FYa(\L, \G_{\min})$ admits a quadratic Gröbner basis and is therefore Koszul.
\end{theo}

Stanley \cite{Stanley_1972} proved that the geometric lattices associated to chordal graphs (i.e. graphs such that every cycle has a chord) are supersolvable. Alternatively, one can associate to a graph $G$ a built lattice $(\L_G , \G_G)$ where $\L_G$ is the lattice associated to $G$ and $\G_G$ is the building set of connected closed subgraphs of $G$. Stanley's original argument also shows that $(\L_G, \G_G)$ is a supersolvable built lattice.  This implies by Theorem~\ref{theomainintro} that its  Feichtner--Yuzvinsky algebra admits a quadratic Gröbner basis. Since the complete graphs are chordal we recover the result of Dotsenko. \\

In \cite{LM_2000}, Losev and Manin introduced moduli spaces of stable curves with marked points of two types, where the points of the first type are not allowed to coincide with any other points, and those of second type are allowed to coincide between them. Those moduli spaces form the components of an object called the ``extended modular operad'', introduced by Losev and Manin in the sequel \cite{LM_2004}. In \cite{Manin_2004}, Manin asked if the cohomology algebras of those moduli spaces are Koszul. By considering the family of chordal graphs $G_{m,n}$, where $G_{m,n}$ has $m+n$ vertices, the first $m$ vertices are neighbors of every vertices and the last $n$ vertices are neighbors only of the first $m$ vertices, one obtains the following result. 
\begin{theo}\label{theoextintro}
The cohomology algebras of the components of the extended modular operad in genus $0$ have quadratic Gröbner bases and are therefore Koszul. 
\end{theo}

In Section \ref{secpreli} we introduce the combinatorial objects needed to understand the rest of the article and we recall some of their known properties.\\

In Section \ref{secmainres} we prove Theorem \ref{theomainintro} and we deduce Theorem \ref{theogminintro}. \\

In Section \ref{secgraph} we concentrate our attention toward supersolvable built lattices associated to chordal graphs, which leads to Theorem \ref{theoextintro}. \\

Finally, in Section \ref{secfurther} we take a step back and give some general comments for further research. \\

\textbf{Acknowledgements.} The author would like to thank Vladimir Dotsenko for his time spent on re-readings and his highly salutary redactional comments. We would also like to express our gratitude toward Russ Woodroofe for his comments on the use of supersolvability.  \\

This research is part of the author's PhD. The author would like to thank Eva Maria Feichtner and Karim Adiprasito for agreeing to be his thesis referees. \\

This research was supported by the University of Strasbourg Institute for Advanced Study through the French national program ``Investment for the future'' [IdEx-Unistra, fellowship USIAS-2021-061 of V.~Dotsenko] and by the French national research agency [grant ANR-20-CE40-0016].

\section{Preliminaries}\label{secpreli}

In this section we introduce the main combinatorial objects which will be used throughout this paper. 
\subsection{Geometric lattices and matroids}\label{subsecmat}
\begin{madef}[Lattice]
A finite poset $\L$ is called a \textit{lattice} if every pair of elements in $\L$ admits a supremum and an infimum.
\end{madef}
The supremum of two elements $G_1, G_2$ is denoted by $G_1 \vee G_2$ and called their \textit{join}, while their infimum is denoted by $G_1 \wedge G_2$ and called their \textit{meet}.
\begin{rmq}
Since $\L$ is supposed to be finite, having supremums and infimums for pairs of elements implies having supremums and infimums for any subset $S$ of $\L$, which will be denoted by $\bigvee S$ and $\bigwedge S$ respectively. As a consequence, every lattice admits an upper bound (the supremum of $S = \L$) and a lower bound (the infimum of $S = \L$) which will be denoted by $\un$ and $\zero$ respectively.
\end{rmq}
\begin{madef}[Geometric lattice]\label{defgeolat}
A finite lattice $(\L, \leq)$ is said to be \textit{geometric} if it satisfies the following properties:
\begin{itemize}
\item For every pair of elements $G_1 \leq G_2$, all the maximal chains of elements between $G_1$ and $G_2$ have the same cardinal. \textit{(Jordan-Hölder property)}
\item The rank function $\rho: \L \rightarrow \N$ which assigns to any element $G$ of $\L$ the cardinal of any maximal chain of elements from $\zero$ to $G$ (not counting $\zero$) satisfies the inequality
\begin{equation*}
\rho(G_1 \wedge G_2) + \rho(G_1 \vee G_2) \leq \rho(G_1) + \rho(G_2)
\end{equation*}
for every $G_1$, $G_2$ in $\L$. \textit{(Sub-modularity)}
\item Every element in $\L$ can be obtained as the supremum of some set of atoms (i.e. elements of rank~$1$). \textit{(Atomicity)}
\end{itemize}
\end{madef}
For any geometric lattice $\L$ we will denote by $\At(\L)$ its set of atoms. For any element $F$ in $\L$ we will denote by $\At_{\leq}(F)$ the set of atoms of $\L$ which are below $F$. \\

One of the  reasons to study this particular class of lattices is that the intersection poset of any hyperplane arrangement is a geometric lattice. In fact, one may think of geometric lattices as a combinatorial abstraction of hyperplane arrangements. In addition, this object is equivalent to the datum of a loopless simple matroid via the lattice of flats construction and therefore it has connections to many other areas in mathematics (graph theory for instance). \\

Let us describe in more details the correspondence between simple loopless matroids and geometric lattices. There are several equivalent definitions of matroids. We refer to  \cite{welsh_matroid_1976} for more details. 
\begin{madef}[Matroids via independent subsets]
A matroid is a pair of a finite set $E$ and a set $\I$ of subsets of $E$ (the ``independent'' subsets) satisfying the axioms
\begin{itemize}
\item For any $I$ in $\I$, every subset of $I$ belongs to $\I$.
\item For any $I$, $J$ in $\I$, if $\#J > \#I$ there exists an element $a$ in $J$ and not in $I$ such that $I \cup \{a\}$ is independent.
\end{itemize}
\end{madef}
\begin{madef}[Matroids via closure operator]
A matroid is a pair of a finite set $E$ and an application (the ``closure operator'')
\begin{equation*}
\sigma: \P(E) \longrightarrow \P(E)
\end{equation*}
satisfying the axioms 
\begin{itemize}
\item For any $X \in \P(E)$ we have $X \subseteq \sigma(X)$.
\item For any $X \subseteq Y \in \P(E)$ we have $\sigma(X) \subseteq \sigma(Y)$. 
\item For any $X \in \P(E)$ we have $\sigma(\sigma(X)) = \sigma(X)$. 
\item For any $X \in \P(E)$ and $a, b \in E$, if $a \in \sigma(X \cup \{b\}) \setminus \sigma(X)$ then $b \in \sigma(X \cup \{a\}) \setminus \sigma(X)$. 
\end{itemize}
\end{madef}
\begin{madef}[Matroids via circuits]\label{defmatcirc}
A matroid is a pair of a finite set $E$ and a set $\Circ$ of subsets of $E$ (the ``circuits'') satisfying the axioms
\begin{itemize}
\item The empty set is not a circuit. 
\item If $C_1 \subseteq C_2 \in \Circ$ then $C_1 = C_2$. 
\item If $C_1, C_2 \in \Circ$, $C_1 \neq C_2$ and $e \in C_1 \cap C_2$ then there exists a circuit $C \subseteq C_1 \cup C_2\setminus \{e\}$.
\end{itemize}
\end{madef}
One can replace the last axiom by a stronger version which we will use later in this article. 
\begin{multline}\label{axcircuit}
\textrm{If } C_1, C_2 \in \Circ, \, e \in C_1 \cap C_2, \, f \in C_1 \setminus C_2, \\ \textrm{ then there exists a circuit } C \subseteq C_1 \cup C_2 \setminus \{e\} \textrm{ containing } f. 
\end{multline}
One passes from the independent subset definition to the circuit definition by defining a circuit as a minimal dependent subset. One passes from the circuit definition to the closure definition by putting 
\begin{equation}\label{eqclosurecircuit}
\sigma(X) \coloneqq X \cup \{x  \, | \, \exists C \in \Circ \textrm{ with } C \subseteq  X \cup \{x\} \textrm{ and } x\in C\}. 
\end{equation}
A matroid $(E, \I)$ is said to be simple loopless if every subset of $E$ of cardinal less than two is independent. A flat of a matroid $M = (E, \sigma)$ is a subset $F \subseteq E$ such that $\sigma(F)$ is equal to $F$. The set of flats of $M$ denoted by $\L_M$ ordered by inclusion is a geometric lattice with meet given by the intersection. Conversely if $\L$ is a geometric lattice then the datum $(E, \sigma)$ where $E$ is the set of atoms of $\L$ and $\sigma$ is the map defined by 
\begin{equation*}
\sigma(X) = \bigcap_{\substack{F \in \L \\ X \subset \At_{\leq}(F)}} \At_{\leq}(F) 
\end{equation*}
is a simple loopless matroid. Those two constructions are inverse to each other on simple loopless matroids. In the sequel we will freely identify an element of some geometric lattice with the set of atoms below this element. For instance if $G_1$ and $G_2$ are two elements of some geometric lattice $\L$ then $G_1 \cup G_2$ will mean $\At_{\leq}(G_1) \cup \At_{\leq}(G_2)$. Finally, notice that by definition, for any subset $S \subseteq \L$ with $\L$ some geometric lattice we have 
\begin{equation}\label{eqclosurejoin}
\bigvee S = \sigma(\bigcup_{X \in S} X) 
\end{equation}
where $\sigma$ is the closure operator of the associated matroid. \\

Here is a list of some important well-known geometric lattices.
\begin{ex}\label{exgeolatt}
\leavevmode
\begin{itemize}
\item If $X$ is any finite set, the set $\mathcal{P}(X)$ of subsets of $X$ ordered by inclusion is a geometric lattice with join the union and meet the intersection. It is the intersection lattice of the hyperplane arrangement of coordinate hyperplanes in $\Cx^{X}$. Those geometric lattices are called boolean lattices and denoted by $\mathcal{B}_{X}$.
\item If $X$ is any finite set, the set $\Pi_{X}$ of partitions of $X$ ordered by refinement is a geometric lattice. It is the intersection lattice of the so-called \textit{braid} arrangement which consists of the diagonal hyperplanes $\{z_i = z_j\}$ in $\Cx^{X}$. Those geometric lattices are called partition lattices.
\item If $G = (V, E)$ is any graph one can construct the graphical matroid $M_{G}$ associated to $G$ and then consider $\L_G$ the lattice of flats associated to $M_{G}$ (see \cite{welsh_matroid_1976} for the details of this construction). Those lattices are said to be \emph{graphical}. This family of geometric lattices contains the two previous ones because $\mathcal{B}_{X}$ is the lattice associated to any tree with edges $X$ and $\Pi_{X}$ is the lattice associated to the complete graph with vertices $X$. For any graph $G = (V, E)$ the geometric lattice $\L_{G}$ is the intersection lattice of the hyperplane arrangement $\{\{z_{u} = z_{v}\}, (u,v) \in E\}$ in $\mathbb{C}^{V}$. 
\end{itemize}
\end{ex}
A fundamental fact about geometric lattices is that an interval of a geometric lattice is a geometric lattice (see \cite{welsh_matroid_1976}). \\

In the rest of this article every lattice will be assumed to be geometric.
\subsection{Building sets and nested sets}
The following definition is due to De Concini-Procesi \cite{de_concini_wonderful_1995} in the realizable case, and Feichtner-Yuzvinsky \cite{FY_2004} in general. 
\begin{madef}[Building set]\label{defbuilding}
Let $\L$ be a geometric lattice. A \textit{building set} $\G$ of $\L$ is a subset of $\L \setminus \{\zero\}$ such that for every element $X$ of $\L$ the morphism of posets
\begin{equation}\label{isobuilding}
\prod_{G \in \max \G_{\leq X}}[\zero, G] \xrightarrow{\vee} [\zero, X]
\end{equation}
is an isomorphism (where $\max \G_{\leq X}$ is the set of maximal elements of $\G \cap [\zero,X]$).
\end{madef}
The elements of $\max \G_{\leq X}$ are called the \emph{factors} of $X$ in $\G$. 

\begin{madef}[Built lattice]
The datum of a lattice $\L$ and a building set $\G$ of $\L$ will be called a \textit{built lattice}. If $\G$ contains $\un$ we say that $(\L, \G)$ is \textit{irreducible}.
\end{madef}
The definition of a building set makes sense for a larger class of posets, as shown in $\cite{FY_2004}$, but in this paper we will restrict ourselves to the case of geometric lattices. In this particular context, building sets are geometrically motivated by the construction of wonderful compactifications for hyperplane arrangement complements. \\

To put it in a nutshell, building sets are sets of intersections of a hyperplane arrangement that one can successively blow up in order to obtain a wonderful compactification of its complement (see \cite{de_concini_wonderful_1995} for more details). Each blowup creates a new exceptional divisor, so the wonderful compactification is equipped with a family of irreducible divisors indexed by $\G$. This family of divisors forms a normal crossing divisor when $\G$ is a building set. \\

There are a few key examples to keep in mind throughout this story.
\begin{ex}\label{exbs}
\leavevmode
\begin{itemize}
\item Trivially, every lattice $\L$ admits $\L\setminus\{\zero\}$ as a building set.
\item Less trivially, every lattice $\L$ also admits a unique minimal building set which consists of all the elements $G$ of $\L$ such that $[\zero,G]$ is not a product of proper subposets.
\item From the definition one can see that a building set of some lattice $\L$ must contain all the atoms of $\L$. If $\L$ is a boolean lattice (see Example \ref{exgeolatt}) then its set of atoms is in fact a building set (the minimal one). This fact characterizes boolean lattices.
\item If $\L$ is the lattice of partitions of some finite set (see Example \ref{exgeolatt}) then the subset of partitions with only one block having more than two elements is a building set of $\L$. This is the minimal building set of $\L$.
\item If $\L$ is a graphical lattice (see Example \ref{exgeolatt}) then the set of elements of $\L$ corresponding to sets of edges which are connected forms a building set of $\L$. We will denote this built lattice associated to a graph $G$ by $(\L_G, \G_G)$. This family of examples contains the two previous ones (by considering totally disconnected graphs for the former and complete graphs for the latter).
\item Alternatively, if $G = (V,E)$ is a graph one can consider the boolean lattice $\B_{V}$. This lattice has a building set made up of the ``tubes'' of $G$, that is sets of vertices of $G$ such that the induced subgraph on those vertices is connected. This leads to the notion of graph associahedra introduced in \cite{CD_2004}.
\end{itemize}
\end{ex}
A key fact about building sets is that any interval $[G_1, G_2]$ in some built lattice $(\L, \G)$ admits an ``induced'' building set which we describe now. We start by introducing a useful notation.
\begin{nota}
For any element $G$ of some lattice $\L$ and a subset $X$ of $\L$, we denote by $G\vee X$ the set of elements of $\L$ which can be obtained as the join of $G$ and some element of $X$.
\end{nota}
\begin{madef}[Induced building set]
Let $G_1 < G_2$ be two elements in some built lattice~$(\L, \G)$. We denote by $\Ind_{[G_1, G_2]}(\G)$ the set $(G_1\vee \G) \cap [G_1, G_2] \setminus \{G_1\}$ and we call it the \emph{induced building set} on $[G_1, G_2]$.
\end{madef}
\begin{lemma}
The subset $\Ind_{[G_1, G_2]}(\G) \subset (G_1, G_2]$ is a building set of $[G_1, G_2]$.
\end{lemma}
We will often write $\Ind(\G)$ instead of  $\Ind_{[G_1, G_2]}(\G)$ if the interval can be deduced from the context.
\begin{proof}
The proof can be found in \cite{Bibby_2021} (Lemma 2.8.5).
\end{proof}
\begin{madef}[Nested set]
Let $(\L, \G)$ be a built lattice. A subset $\S$ of $\G$ is called a \textit{nested set} if for every antichain $\A$ in $\S$ which is not a singleton, the join of the elements of $\A$ does not belong to $\G$. The nested sets of a built lattice $(\L, \G)$ form an abstract simplicial complex denoted $\Ne(\L,\G)$. We denote by $\Ne^{\irr}(\L, \G)$ the set of nested sets of $(\L, \G)$ containing the maximal elements of $\G$. 
\end{madef}
\begin{ex}
A chain of elements in some building set $\G$ is always nested and those are the only nested sets of the maximal building set ($\G = \L \setminus \{\zero \}).$
\end{ex}
Geometrically, nested sets correspond to sets of divisors in the wonderful compactification which have a nontrivial intersection.
\subsection{Supersolvable built lattices}
Recall that a lattice $\L$ is said to be distributive if for every triple $X, Y, Z \in \L$ we have the equality $$X \wedge(Y \vee Z) = (X\wedge Y)\vee (X \wedge Z).$$
A pair of elements $(X,Y)$ in a lattice $\L$ is said to be modular if for every $Z \leq Y$ we have the identity $$Z\vee(X\wedge Y) = (Z\vee X)\wedge Y.$$ 
An element $X$ in a lattice $\L$ is said to be modular if for every $Y$ both the pairs $(X,Y)$ and $(Y,X)$ are modular. \\

The following definition is due to Stanley \cite{Stanley_1972}.
\begin{madef}[Supersolvable lattice]
A lattice $\L$ is said to be \textit{supersolvable} if there exists a maximal chain $M$ of elements of $\L$ such that for every chain $K$ in $\L$ the sublattice generated by $M$ and $K$ is distributive. 
\end{madef}
Stanley proved that for geometric lattices (or more generally semimodular lattices) we have the following equivalence. 
\begin{prop}[Stanley \cite{Stanley_1972}]
A geometric lattice is supersolvable if and only if it has a maximal chain of modular elements. 
\end{prop}
In this article we only consider geometric lattices and we will mostly use the above equivalent characterization. 
\begin{fact}\label{facthereditary}
Supersolvability is a hereditary condition (meaning it is stable by taking intervals) because if $G$ is some element in some supersolvable lattice $\L$ with maximal chain of modular elements $$\zero = M_1 < ... < M_n = \un,$$ the maximal chains $$M_1\wedge G \leq ... \leq M_n \wedge G$$ and $$M_1 \vee G \leq ... \leq M_n \vee G$$ are maximal chains (with possibly multiple occurencies) of modular elements of $[\zero, G]$ and $[G, \un]$ respectively (see \cite{Stanley_1972}).
\end{fact}
We introduce the following variant for built lattices. 

\begin{madef}[Supersolvable built lattices]
A built lattice $(\L, \G)$ is said to be supersolvable if it admits a maximal chain $\zero = G_1 < ... < G_n = \un$ of modular elements in $\G$ such that for any element $G$ in $\G$, the element $G_i \wedge G$ belongs to $\G \cup \{\zero\}$ for all $i\leq n$.
\end{madef}
By Fact \ref{facthereditary}, supersolvability for built lattices is a hereditary condition (meaning it is stable by taking intervals and induced building set). 

\begin{ex}
Let $\B_{4}$ be the boolean lattice of $\{1, 2, 3 ,4\}$. If we put $$\G \coloneqq \{\{1\}, \{2\}, \{3\}, \{4\}, \{1,2\}, \{1,2,3\}, \{1,2,3,4\}, \{2,3,4\}, \{2,3\}\}$$ then the built lattice $(\B_{4}, \G)$ is supersolvable. Indeed the chain $$\{1\} \subset \{1, 2\} \subset \{1, 2 ,3\} \subset \{1, 2, 3, 4\}$$ is a maximal chain of modular elements in $\G$ (all the elements of $\B_4$ are modular), and we have $$\{2,3,4\} \wedge \{1,2,3\} = \{2, 3\} \in \G.$$ On the contrary if one puts $\G' \coloneqq \G \setminus \{2, 3\} \cup \{3,4\}$ then $(\B_4, \G')$ is not a supersolvable built lattice because we have $$\{1,2,3\} \wedge \{2, 3, 4\} = \{2, 3\} \notin \G.$$ and any maximal chain of elements in $\G$ must contain either $\{1, 2 ,3\}$ or $\{2, 3,4\}$. 
\end{ex}

One can immediately see that if $\L$ is a supersolvable lattice then $(\L, \G_{\max})$ is a supersolvable built lattice. In Section \ref{secmainres} and Section \ref{secgraph} we will introduce other large classes of supersolvable built lattices. \\

Let $(\L, \G)$ be a supersolvable built lattice with some chosen maximal chain of modular elements $\omega = \{ \zero = G_1 < ... < G_n = \un\}$. For any $G$ in $\G$ and any $G'$ in $\Ind_{[G, \un]}(\G)$ we denote by $\d_{\omega,G}(G')$ the coatom in the maximal chain of modular elements induced by $\omega$ on $[G, G']$ (see Fact \ref{facthereditary}). An element of the form $\d_{\omega, G}(G')$ will be called an \textit{initial segment} of $G'$ relative to $G$. In practice we will drop $\omega$ from the notation. If $G$ is equal to $\zero$ we also drop it from the notation. In the sequel whenever we introduce a supersolvable built lattice we implicitly choose a particular maximal chain of modular elements of this built lattice. \\

To conclude this subsection we prove a small general lemma which will be useful later on. 
\begin{lemma}\label{lemmamodularcircuit}
Let $\L$ be a geometric lattice,  $G$ a modular element of $\L$ and $C$ a circuit in $\L$. At least one of the following propositions is true:
\begin{itemize}
\item $C \subset G$. 
\item $C \cap G = \emptyset$.
\item There exists a circuit $C'$ such that $C'\cap G$ is a singleton $H$ and $C'\setminus\{H\}$ is included in $C\setminus (G\cap C)$.
\end{itemize}
\end{lemma}
\begin{proof}
Assume the first two propositions are not true. Let us denote 
\begin{align*}
&I \coloneqq C \cap G, \\
&J \coloneqq C \setminus I.
\end{align*}
Let $H$ be any element of $I$, which is not empty by assumption. By modularity of $G$ we have
\begin{align*}
G\wedge (\sigma(I \setminus \{H\}) \vee \sigma(J)) = \sigma(I \setminus\{H\}) \vee (G\wedge \sigma(J)).  
\end{align*}
Since $J$ is not empty $I$ is independent and the element on the left has rank at least $\#I$. Since $\sigma(I \setminus\{H\})$ only has rank $\#I -1$ we must have $$G\wedge \sigma(J) \neq \zero.$$ By Formula \eqref{eqclosurecircuit} this implies that we have the desired circuit $C'$. 
\end{proof}
\subsection{The Feichtner--Yuzvinsky rings}
\begin{madef}
For every built lattice $(\L, \G)$ we define the Feichtner--Yuzvinsky graded commutative ring $\FYa(\L, \G)$ by
\begin{equation*}
\FYa(\L, \G) = \Z[x_G, \,G \in \G]/ \I_{\aff},
\end{equation*}
with all the generators in degree $2$, and $\I_{\aff}$ the ideal generated by elements
\begin{equation*}
\sum_{G \geq H}x_G
\end{equation*}
for every atom $H$, and elements
\begin{equation*}
\prod_{G \in X}x_G
\end{equation*}
for every set $X \subset \G$ which is not nested.
\end{madef}
In the realizable case, the ring $\FYa(\L, \G)$ is the cohomology ring of the wonderful compactification associated to the building set $\G$ (see \cite{de_concini_wonderful_1995} for the computation of the cohomology ring). Those rings were generalized to arbitrary built lattices by Feichtner and Yuzvinsky in \cite{FY_2004}. \\

The Feichtner--Yuzvinsky rings admit another important presentation.
\begin{prop}
For every built lattice $(\L, \G)$ we have the presentation
\begin{equation*}
\FYa(\L, \G) \simeq \Z[h_G, \, G \in \G] / \I_{\wond}
\end{equation*}
where $\I_{\wond}$ is the ideal generated by relations
\begin{equation*}
h_H
\end{equation*}
for every atom $H$ and
\begin{equation*}
\prod_{G' \in \A} (h_G - h_{G'})
\end{equation*}
for every $G \in \G$ and $\A$ an antichain in $\G$ such that $\bigvee \A$ is equal to $G$. The change of variable between the last presentation and the defining presentation is given by $$h_G = \sum_{G' \geq G}x_{G'}.$$
\end{prop}
This presentation appeared first in \cite{EHKR_2010} for the braid arrangement and in \cite{Backman_Spencer_Eur_2020} for general maximal building sets. It is widely used in \cite{Pagaria_2021}. In this article we will use exclusively this presentation. 
\begin{proof}
The proof can be found in \cite{Pagaria_2021} (Theorem 2.9).
\end{proof}
In \cite{FY_2004}, the authors address the issue of finding a Gröbner basis for $\FYa_{\aff}(\L, \G)$ (see \cite{BW_1993} for a reference on Gröbner bases) and they show that when considering any linear order on generators refining the reverse order on $\G$, although the elements  defining $\I_{\aff}$ do not form a Gröbner basis in general, one can still describe a fairly manageable Gröbner basis.
\begin{theo}[Feichtner--Yuzvinsky]\label{theogrobnerFY}
Elements of the form $(\prod_{G \in \S}x_G) h_{G'}^{\rho(G') - \rho(\bigvee S)}$ with $\S$ any nested set and $G'$ any element of $\G$ satisfying $G' > \bigvee \S$ ,  together with the usual $\prod_{G \in X} x_G$ for every non-nested set $X$, form a Gröbner basis of $\FYa_{\aff}(\L, \G)$ for any linear order on generators refining the reversed order of $\L$. The normal monomials with respect to this Gröbner basis are monomials of the form
\begin{equation*}
x_{G_1}^{\alpha_1}... x_{G_n}^{\alpha_n}
\end{equation*}
where the $G_i$'s form a nested set $\S$ and for every $i \leq n$ we have $\alpha_i < \rk [\bigvee \S_{< G_i}, G_i]$.
\end{theo}
\begin{proof}
The proof can be found in \cite{FY_2004}.
\end{proof}
Those Gröbner bases are almost never quadratic. The rest of the paper will be devoted to proving  the following theorem.
\begin{theo}\label{theomain}
Let $(\L, \G)$ be a supersolvable built lattice. The Feichtner--Yuzvinsky ring $\FYa(\L, \G)$ admits a quadratic Gröbner basis and is therefore Koszul. 
\end{theo}
The proof of this theorem will be carried out in Subsections \ref{secorder} and \ref{secbij}. The following proposition is a small step toward Theorem \ref{theomain} which we will need later on. 
\begin{prop}\label{propquad}
Let $(\L, \G)$ be a supersolvable built lattice. The Feichtner--Yuzvinsky ring $\FYa(\L, \G)$ is quadratic.
\end{prop}
\begin{proof}
It is enough to prove that the nested set complex of $(\L, \G)$ is flag, meaning that for any anti-chain $G_1, ... , G_n$ in $\G$ with $n\geq 2$, if $G_1 \vee ... \vee G_n$ belongs to $\G$ then there exists $i \neq j \leq n$ such that $G_i \vee G_j$ belongs to $\G$. Assume the contrary is true and there exist an anti-chain $G_1, ..., G_n$ such that we have 
\begin{equation*}
G_1 \vee ... \vee G_n  \in \G 
\end{equation*}
\begin{center}
and
\end{center}
\begin{equation*}
G_i \vee G_j \notin \G \, \, \forall i \neq j \leq n.
\end{equation*}
By restricting to a smaller interval we can assume $G_1 \vee ... \vee G_n  = \un$. By Fact \ref{facthereditary}, the element $\d(\un)\wedge (G_i \vee G_j)$ is either equal to $G_i \vee G_j$ or is covered by $G_i \vee G_j$. In the latter case, using the building set isomorphism \eqref{isobuilding} we see that either $G_i$ or $G_j$ is below $\d(\un)$. As a consequence we see that there is at most one integer $i \leq n$ such that $G_i$ is not below $\d(\un)$ (in fact there is exactly one such $i$). By reordering let us assume that this integer is $n$. By atomicity there exist an element $X < G_n$ such that we have $G_1 \vee ... \vee G_{n-1} \vee X = \d(\un)$. If $G'_1, ... G'_k$ are the factors of $X$ in $\G$, the anti-chain $G_1, ... , G_k, G'_1, ..., G'_k$ is a new counter-example to the flagness of the nested set complex of $(\L,\G)$. We get a contradiction by reiterating this process. 
\end{proof}
\subsection{The operadic structure on built lattices}\label{subsecop}
In order to prove Theorem \ref{theomain} we will use the operadic structure on built lattices introduced in \cite{Coron_2022}. Let us quickly summarize this construction, referring to the latter article for more details. From now on, all nested sets are supposed to contain the maximal elements of the building set they live in. \\

For any built lattice $(\L, \G)$ and any nested set $\S$ in $\G$ we have maps of algebras
\begin{equation}\label{eqmorstr}
\begin{array}{ccc}
\FYa(\L, \G) & \xrightarrow{\FYa(\S)} &\bigotimes_{\G \in \S}\FYa([\bigvee \S_{<G}, G],\Ind(\G)) \\
h_{G'} & \longrightarrow & 1\otimes ... \otimes 1 \otimes h_{(\bigvee S_{<G})\vee G'} \otimes 1 \otimes ... \otimes 1 
\end{array} 
\end{equation}
(where $G$ is the minimal element of $\S$ such that we have $G' < G$). In \cite{Coron_2022} we show that the collection of Feichtner--Yuzvinsky algebras $\{\FYa(\L, \G), (\L, \G)\}$ together with the above morphisms can be formalized as an operad over a certain Feynman category (Kaufmann-Ward~\cite{kaufmann_feynman_2017}) denoted $\LBS$ having as objects the built lattices and morphisms the nested sets. The key ingredient is the definition of an associative composition of nested sets : for every nested set $\S \in \Ne^{\irr}(\L, \G)$ and every collection of nested sets $$(\S_G, \in \Ne^{\irr}([\bigvee \S_{<G}, \un ], \Ind(\G)))_{G \in \S}$$ one can define $\S\circ(\S_{G})_{G} \in \Ne^{\irr}(\L, \G)$ such that the operation $\circ$ is associative (this is the composition of morphisms in $\LBS$). In the case of the maximal building set, the composition of nested sets is just the concatenation of chains. \\

The operad $\FY^{\vee} = \{\FYa^{\vee}(\L, \G), (\L, \G)\}$ with structural morphisms given by the linear dual of morphisms \eqref{eqmorstr} admits a quadratic presentation with one generator of top degree (the degree map) in each arity (each built lattice $(\L, \G)$). The ``monomials'' in $\FY^{\vee}$ are all possible operadic products of those generators. Since we have only one generator in each arity, the monomials living in some Feichtner--Yuzvinsky ring $\FYa(\L, \G)$ are in bijection with the nested sets of $(\L, \G)$. In \cite{Coron_2022} we construct a Gröbner basis machinery for operads over $\LBS$, via the introduction of a notion of shuffle $\LBS$-operads, governed by another Feynman category having as objects the directed built lattices (built lattices together with an order on the atoms of $\L$). One can compute an (operadic) Gröbner basis for $\FY^{\vee}$ and describe the associated (operadic) normal monomials as follow.  Let $(\L, \G, \vartriangleleft)$ be a directed built lattice. The total order on atoms $\vartriangleleft$ defines an EL-labelling $\EL$ by putting
\begin{equation*}
\EL(X\prec Y ) = \min \,\{H \in \At(\L) \, | \, X \vee H = Y\}.
\end{equation*}
We refer to \cite{wachs_poset_2006} for a reference on EL-labellings. For any two elements $X< Y$ in $\L$ and $k$ some positive integer less than $\rk(Y) - \rk(X)$ we define $\omega^k_{X, Y, \EL}$ to be the truncation at height $k$ of the unique maximal chain with increasing $\EL$-labels between $X$ and $Y$. More precisely, if this unique maximal chain is
\begin{equation*}
X_0 = X \prec X_1 \prec ... \prec X_n = Y,
\end{equation*}
then $\omega^k_{X, Y, \EL}$ is the chain
\begin{equation*}
X_0 \prec ... \prec X_k.
\end{equation*}
If there is no ambiguity on the  EL-labelling we will drop it from the notation. If the maximal chain is not truncated, i.e. $k ~=~\rk(Y) -~\rk(X)$, we simply omit the superscript. In addition, to any maximal chain $\omega = \zero \prec X_1 \prec ... \prec X_n $ in $\L$ one can associate a nested set
\begin{equation*}
\S(\omega) = \{X_1\} \circ ... \circ \{X_n\}.
\end{equation*}
Notice that this is well defined even if the $X_i$'s do not belong to the building set $\G$ because every $X_i$ is an atom in $[X_{i-1}, \un]$ and therefore must belong to the induced building set on the interval $[X_{i-1}, \un]$. Finally, in \cite{Coron_2022} we show that the operadic normal monomials are represented by the nested sets of the form
\begin{equation*}
\S = \S'\circ (\S(\omega^{k_G}_{\tau_{\S}(G), G}))
\end{equation*}
where we have $k_G < \rk([ \tau_{\S'}(G) , G ]) - 1$ for all $G$ in $\S'$, except for $k_{\un}$ which is less or equal than $\rk([\tau_{\S'}(\un), \un])$.
\section{The main results}\label{secmainres}
This section is devoted to the proof of Theorem \ref{theomain}. In the first subsection we define an order on the generators of the Feichtner--Yuzvinsky algebras of supersolvable lattices and we compute the normal monomials of weight $2$ associated to this order. In the next subsection we define a bijection between the algebraic normal monomials associated to the latter order and the operadic normal monomials introduced in Subsection \ref{subsecop}. The construction of this bijection will be done by induction and using the operadic structure. By a dimension argument this bijection will show that the algebraic normal monomials form a basis of the Feichtner--Yuzvinksky algebras, which will show that the set of weight $2$ relations forms a Gröbner basis of the Feichtner--Yuzvinsky algebras. 
\subsection{The order on generators and the normal monomials of weight~2}\label{secorder}
\begin{madefprop}\label{deforder}
Let $(\L, \G)$ be a supersolvable built lattice. The transitive closure of the relations 
\begin{equation}\label{eqorder}
\d^{k}(G) \dashv G \dashv G'
\end{equation}
for all $k$ and all $G,G' \in \G$ with $G' \leq G$ and $G'$ not an initial segment of $G$, is anti-symmetric and thus defines a partial order. 
\end{madefprop}
\begin{proof}
One can define an explicit total order containing the relations \eqref{eqorder} as follow. Let $\vartriangleleft$ be a total order on the atoms of $\L$ extending the relations $H\dashv H'$ for all pairs of atoms $H, H'$ such that there exists an integer $k$ satisfying 
\begin{equation*}
H \leq \d^{k}(G) \textrm{ and } H \nleq \d^{k}(G).
\end{equation*}
For any element $G$ in $\G$ let us denote by $\w(G)$ the word with letters $\At_{\leq}(G)$ written in increasing order. 
We define a total order on $\G$, also denoted $\vartriangleleft$, by putting
\begin{equation*}
G \vartriangleleft G' \Leftrightarrow \w(G) \textrm{ is less than } \w(G') \textrm{ for the lexicographic order.} 
\end{equation*}
This order contains relations \eqref{eqorder}. 
\end{proof}
In the sequel whenever we introduce a supersolvable built lattice we implicitly choose an associated total order on $\G$ as in the above proof. 
\begin{prop}\label{propnormal2}
Let $(\L, \G)$ be a supersolvable built lattice and let $\alpha = h_{G_1}h_{G_2}$ be a monomial in $\FYa(\L, \G)$ with $G_1 \vartriangleleft G_2$. Let us denote by $G$ the join $G_1 \vee G_2$. The monomial $\alpha$ is normal if and only if one of the three following conditions is verified.
\begin{itemize}
\item The element $G$ does not belong to $\G$.
\item The element $G$ belongs to $\G$ and $G_1$ is an initial segment of $G_2 = G$.
\item The element $G$ belongs to $\G$, $G_2 \ngeq G_1$, $G_2$ is not covered by $G$ and we have $G_1 = \d^{k}(G)$ where $k$ is the maximal integer satisfying
\begin{equation*}
\d^{k}(G) \vee G_2 = G.
\end{equation*}
\end{itemize}
\end{prop}

\begin{proof}
We have an obvious bijection between the monomials described in the above proposition and the normal monomials given by \ref{theogrobnerFY}, sending $h_{\d^{i_{G'}}(G)}h_{G'}$ to $x_{G'}x_{G}$ if $G'$ is not covered by $G$, sending $h_{\d(G)}h_G$ to $x_G^2$ and sending $h_{G_1}h_{G_2}$ to $x_{G_1}x_{G_2}$ when $\{G_1, G_2\}$ is a nested set. By a dimension argument it is enough to prove that the normal monomials of weight $2$ with respect to $\vartriangleleft$ are included in the monomials described in the proposition. \\

Assume $G = G_1 \vee G_2$ belongs to $\G$. If $G \vartriangleleft G_1$ and $G \vartriangleleft G_2$ then $h_{G_1}h_{G_2}$ is the leading term of the relation
\begin{equation*}
(h_{G}-h_{G_1})(h_{G}- h_{G_2}).
\end{equation*}
If $G_1 \totoreq G \vartriangleleft G_2$ and $G_2$ is covered by $G$ then $h_{G_1}h_{G_2}$ is the leading term of the relation 
\begin{equation*}
(h_{G} - h_{G_1})(h_{G} - h_{G_2}) - h_{G}(h_{G} - h_{G_2}). 
\end{equation*}
Finally, if $G_1 \totoreq G \vartriangleleft G_2$ then $G_1$ is some initial segment $\d^{k}(G)$. If $k$ is not the maximal integer such that we have $$\d^k(G)\vee G_2 = G,$$
we see that $h_{G_1}h_{G_2}$ is the leading term of the relation
\begin{equation*}
(h_{G}-h_{G_1})(h_{G}-h_{G_2}) - (h_{G}-h_{\d^{k + 1}(G)})(h_{G} - h_{G_2}).
\end{equation*}

\end{proof}
Notice that the normal monomials do not really depend on $\vartriangleleft$ but only on the chosen maximal chain of modular elements. We next come to an important lemma which highlights a first connection between our normal monomials and supersolvability. 
\begin{lemma}\label{lemmasupernormal}
Let $(\L, \G)$ be a supersolvable built lattice. Let $G_1$ and $G_2$ be two non-comparable elements of $\G$ such that we have $G_1 \vartriangleleft G_2$, $G_1 \vee G_2 \in \G$, $G_2$ is not covered by $G_1\vee G_2$ and $G_1$ is an initial segment of $G_1 \vee G_2$. Then $h_{G_1}h_{G_2}$ is a normal monomial if and only if $G_1 \wedge G_2 \leq d(G_1)$. 
\end{lemma}
The forward statement is always true but for the converse we need the supersolvability hypothesis. For instance consider $\L$ the graphical lattice associated to a $5$-cycle and number the edges (i.e. the atoms) from $1$ to $5$. If we pick $G_1 = \{1, 2, 3 \}$ and $G_2 = \{4, 5\}$ in the maximal building set, then we have $G_1 \wedge G_2 = \zero$ but $h_{G_1}h_{G_2}$ is not normal because we have $\d(G_1) \vee G_2 = G_1 \vee G_2 = \un$. 
\begin{proof}
With the hypothesis on $G_1$ and $G_2$ we have 
\begin{align*}
h_{G_1}h_{G_2} \textrm{ is normal } &\Leftrightarrow \d(G_1)\vee G_2 < G_1\vee G_2 \\
& \Leftrightarrow G_1 \nleq \d(G_1) \vee G_2 \\
& \Leftrightarrow G_1 \wedge (\d(G_1) \vee G_2) < G_1 \\
&\Leftrightarrow \d(G_1) \vee (G_1 \wedge G_2) < G_1  \\
&\Leftrightarrow G_1 \wedge G_2 \leq \d(G_1).
\end{align*} 
The fourth equivalence comes from the fact that by supersolvability $G_1$ is modular in the interval $[\zero, G_1 \vee G_2]$.
\end{proof}

\subsection{A bijection between algebraic normal monomials and operadic normal monomials}\label{secbij}
Let $(\L, \G)$ be a supersolvable built lattice and let $\vartriangleleft$ be a total order on the atoms of $\L$ obtained as in the proof of Definition/Proposition \ref{deforder}. Let us denote $\ANM(\L, \G, \vartriangleleft)$ the algebraic normal monomials with respect to the order $\vartriangleleft$ and the relations of weight $2$, i.e. is the set of monomials in $\FYa(\L, \G)$ which are not divisible by the leading term of some relation of weight $2$. Let us denote by $\ONM(\L, \G, \vartriangleleft)$ the set of operadic normal monomials with respect to the order on atoms $\vartriangleleft$ (see Subsection \ref{subsecop}). Notice that in the supersolvable case the nested set  $\S(\omega^{k_G}_{\tau_{\S}(G), G})$ is simply some truncation of the maximal chain $\{G > \d_{\tau_{\S}(G)}(G) > \d^2_{\tau_{\S}(G)}(G) > ...  > \tau_{\S}(G)\}$. This means that the operadic normal monomials do not depend on the particular choice of $\vartriangleleft$ but only on the choice of the maximal chain of modular elements of $(\L, \G)$. Such choice being implicit we will drop $\vartriangleleft$ from all notations.

\subsubsection{From operadic normal monomials to algebraic normal monomials}
We first define maps
\begin{equation*}
\ONM (\L, \G) \xrightarrow{\Phi_{\L, \G}} \ANM(\L, \G).
\end{equation*}
by induction on the rank of $\L$. Let $(\L, \G)$ be a supersolvable built lattice. For any element $G \neq \un \in \G$ and any algebraic monomial $\alpha  = \prod_{G'\in I} h_{G'}$ in $\FYa([G, \un], \Ind(\G))$ we define the algebraic monomial in $\FYa(\L, \G)$:

\begin{equation*}
\Supp_G(\alpha) = \prod_{\substack{G' \in I \\ G' \in \G}}h_{\d^{i_{G', G}}(G')} \prod_{\substack{G' \in I \\ G' \notin \G}} h_{G'^{\bot}}
\end{equation*}
where for any $G' \in \G$, $i_{G', G}$ is the biggest integer such that we have $\d^{i_{G', G}}(G')\vee G = G'$, and for any $G' \notin \G$, the element $G'^{\bot}$ is the factor of $G'$ in $\G$ different from $G$. Finally, for any operadic normal monomial $\S$ in some supersolvable built lattice $(\L, \G, \vartriangleleft)$ with $G$ the maximal element of $\S'$ for the order $\vartriangleleft$ we define by induction (on both the cardinal of $\S$ and the rank of $\L$) the map
\begin{equation*}
\Phi_{(\L, \G)}(\S)  \coloneqq  \Supp_{G}(\Phi_{[\G, \un], \Ind(\G)}(G\vee \S_{\nleq G})) \Phi_{[\zero, G], \Ind(\G)}(\S_{\leq G}).
\end{equation*}
initialized on empty nested sets by
\begin{equation*}
\Phi(\emptyset) = h_{\un} h_{\d(\un)} ... h_{\d^{\rk \L - 2}(\un)}.  
\end{equation*}
One can check that $G\vee \S_{\nleq G}$ is an operadic monomial so our map is well-defined. This map sends an operadic normal monomial to some algebraic monomial, which will turn out to be normal but we will not need this fact.

\subsubsection{From algebraic normal monomials to operadic normal monomials}
We are concerned with finding an inverse for $\Phi$. Let us define a candidate
\begin{equation*}
\ANM(\L, \G) \xrightarrow{\Psi_{\L, \G}} \ONM(\L, \G)
\end{equation*}
by induction on the rank of $\L$ and the weight of the monomial. We will drop the built lattice from the notation if it can be deduced from the context. We initialize with 
\begin{equation*}
\Psi(1) = \omega_{\zero, \un} = \{ \un > \d(\un) > ... > \zero \}.
\end{equation*}
Let $\alpha = \prod_{i} h_{G_i}$ be some algebraic normal monomial in $\FYa(\L, \G)$ and let us denote by $G$ the maximum of the $G_i's$ with respect to $\vartriangleleft$. If $G = \un$ then by Proposition~\ref{propnormal2} we see that all the $G_i$'s except $G$ are below $\d(\un)$ so we put
\begin{equation*}
\Psi_{\L, \G}(\alpha) = \Psi_{[\zero, \d(\un)], \Ind(\G)}(\alpha/h_{\un}),
\end{equation*}
where $ \Psi_{[\zero, \d(\un)], \Ind(\G)}(\alpha/h_{\un})$ is viewed as a normal monomial in $(\L, \G)$.
If $G$ is different from~$\un$, let us denote respectively
\begin{align*}
&\alpha_{\leq G} = \prod_{G_i \leq G}  h_{G_i}, \\
&\alpha_{\nleq G} = \prod_{G_i \nleq G} h_{G_i}.
\end{align*}
We put
\begin{equation}\label{eqmondef}
\Psi_{\L, \G}(\alpha) \coloneqq \{G\}\circ (\Psi_{[\zero, G], \Ind(\G)}(\alpha_{\leq G}), \Psi_{[G, \un], \Ind(\G)}(G \vee \alpha_{\nleq G})),
\end{equation}
with $G \vee \alpha_{\nleq G}$ a notation for $\prod_{G_i \nleq G}h_{G\vee G_i}$. One can check that this defines an operadic normal monomial. \\

As a side remark let us remind the reader that we have 
\begin{equation*}
\FYa(\{G\}) (\alpha) \coloneqq \alpha_{\leq G}  \otimes G \vee \alpha_{\nleq G},
\end{equation*}
so we are in fact using again the (co)operadic structure on the Feichtner--Yuzvinsky rings. \\

We must prove that the monomial $G\vee \alpha_{\nleq G}$ is normal in $\FYa([G, \un], \Ind(\G))$. This is implied by the following lemma. 
\begin{lemma}\label{lemmamain}
Let $(\L, \G)$ be a supersolvable built lattice and let $G_1 \vartriangleleft G_2 \vartriangleleft G_3$ be elements in $\G$. If $h_{G_1}h_{G_2}$, $h_{G_1}h_{G_3}$ and $h_{G_2}h_{G_3}$ are normal then $h_{G_1 \vee G_3}h_{G_2 \vee G_3}$ is normal in $\FYa([G_3, \un], \Ind(\G))$.
\end{lemma}
This is the technical core of the article. The statement is not true in general without the supersolvability condition, as shown by the following example. Let $\L$ be the graphical lattice associated to a $6$-cycle with edges $\{1, ..., 6\}$. Consider the elements $G_1 = \{1, 2\}, G_2 = \{3,4\}$ and $G_3 = \{5,6\}$. One can quickly check that the monomials $h_{G_1}h_{G_2}$, $h_{G_1}h_{G_3}$ and $h_{G_2}h_{G_3}$ are normal in $\FYa(\L, \G_{\max})$. However, $h_{G_1\vee G_3}h_{G_2 \vee G_3}$ is not normal in $\FYa([G_3, \un], \G_{\max})$ for two reasons: $G_2 \vee G_3$ is covered by $G_1 \vee G_2 \vee G_3 = \un$ and we have $$\{1\}\vee G_2 \vee G_3 = G_1 \vee G_2 \vee G_3.$$ If a built lattice has a small rank it can happen that it satisfies Lemma \ref{lemmamain} without being supersolvable (see Subsection \ref{subsecclassify}).

\begin{proof}
The statement is obvious when two of the $G_i$'s are comparable so we can assume that the elements $G_1, G_2, G_3$ are not comparable. We make a disjunction on whether $G_i\vee G_j$ belongs to $\G$ for $i,j \leq 3$. \\ 

\textbf{Case 1.} $G_i \vee G_j \notin \G$ for all $i \neq j \leq 3$. \\ \\
In this case we have $(G_1\vee G_3)\vee (G_2 \vee G_3) \notin \Ind(\G)$. Indeed, by the proof of Proposition \ref{propquad} the element $G_1\vee G_2\vee G_3$ does not belong to $\G$, and if $G_1\vee G_2\vee G_3$ is equal to $G_1 \vee G$ with $G \in \G$ and $G_1, G$ nested we immediately get $G_1 \vee G_2 = G \in \G$ contradicting the initial hypothesis. \\

\textbf{Case 2.} $G_1 \vee G_2 \notin \G, G_1 \vee G_3 \notin \G, G_2 \vee G_3 \in \G$. \\\\
Let us show that $\{G_1 \vee G_3, G_2 \vee G_3\}$ is a nested anti-chain in $\Ind(\G)$ as in the previous case. By contradiction assume that $G_1 \vee G_2 \vee G_3$ belongs to $\G$. By restriction we can assume $G_1 \vee G_2 \vee G_3 = \un$. Since $h_{G_2}h_{G_3}$ is normal and $G_2 \vee G_3 \in \G$ there exists some integer $k$ such that we have $\d^k(\un)\wedge (G_2 \vee G_3) = G_2$. We have 
\begin{align*}
\d^{k}(\un) &= \d^k(\un) \wedge (G_2 \vee G_1 \vee G_3) \\
&= G_2 \vee (\d^k(\un)\wedge (G_1 \vee G_3)) \\
&= G_2 \vee (\d^k(\un)\wedge G_1) \wedge (\d^k(\un) \wedge G_3) \\
&\leq G_2 \vee G_1 \vee G_2 \\
&= G_1 \vee G_2.
\end{align*}
By nested-ness this implies $\d^k(\un) =G_2$ which contradicts $G_1 \vartriangleleft G_2$. If $G_1\vee G_2 \vee G_3$ belongs to $\Ind(\G)$ but not to $\G$ we immediately get a contradiction as in the previous case. \\

\textbf{Case 3.} $G_1 \vee G_2 \notin \G, G_1 \vee G_3 \in \G, G_2 \vee G_3 \notin \G$. \\ \\ 
This is similar to the previous case. \\

\textbf{Case 4.} $G_1 \vee G_2 \notin \G, G_1 \vee G_3 \in \G, G_2 \vee G_3 \in \G$. \\ \\ 
Once again let us show that $\{G_1 \vee G_3, G_2 \vee G_3\}$ is a nested anti-chain in $\Ind(\G)$. By contradiction assume that $G_1 \vee G_2 \vee G_3$ belongs to $\G$. By restriction we can assume that we have  $G_1 \vee G_2 \vee G_3 = \un$. By assumption there exists an integer $k_1$ such that we have $\d^{k_1}(\un) \wedge (G_1\vee G_3) = G_1$ and an integer $k_2$ such that we have $\d^{k_2}(\un) \wedge (G_2 \vee G_3) = G_2$. Let us denote $k \coloneqq \min (k_1, k_2)$. Let us assume $k_1 \geq k_2$, the other case being symmetric. By modularity of $\d^k(\un)$ and definition of $k$ we have 
\begin{align*}
\d^{k}(\un) &= \d^k(\un) \wedge (G_1 \vee G_2 \vee G_3) \\
&= G_1 \vee (\d^k(\un)\wedge (G_2 \vee G_3)) \\
&= G_1 \vee G_2. 
\end{align*}
This contradicts the fact that $G_1 \vee G_2$ does not belong to $\G$. If $G_1 \vee G_2 \vee G_3$ belongs to $\Ind(\G)$ and not to $\G$ we immediately get a contradiction as in the previous cases. \\

\textbf{Case 5.} $G_1 \vee G_2 \in \G$, $G_1 \vee G_3 \notin \G$ and $G_2 \vee G_3 \notin \G$. \\ \\
We can either have $G_1 \vee G_2 \vee G_3 \notin \G$ or the contrary. In the first case the building set isomorphism 
\begin{equation*}
[\zero ,G_1 \vee G_2 \vee G_3 ] \simeq [\zero, G_1\vee G_2 ]\times[\zero, G_3]
\end{equation*}
immediately gives the result. In the second case we can assume $G_1 \vee G_2 \vee G_3  = \un$. Let us prove that $G_1\vee G_3$ is an initial segment of $G_1\vee G_2\vee G_3$ in $[G_3, \un]$. By assumption there exists an integer $k$ such that we have $\d^{k}(\un) \wedge (G_1 \vee G_2) = G_1$. We will prove the equality 
\begin{equation*}
\d^{k}(\un) \vee G_3 = G_1 \vee G_3.
\end{equation*}
We have 
\begin{align*}
\d^{k}(\un) &= \d^{k}(\un)\wedge(G_1 \vee G_2 \vee G_3) \\
&= G_1 \vee (\d^{k}(\un) \wedge(G_2 \vee G_3)) \, (\textrm{by modularity of } \d^{k}(\un)) \\
&\leq G_1 \vee (\d^{k}(\un) \wedge G_2) \vee (\d^k(\un) \wedge G_3) \, ( \textrm{by nested-ness of} \{G_2, G_3\})\\
&\leq G_1 \vee G_1 \vee G_3 \, (\textrm{by definition of } k) \\
& = G_1\vee G_3. 
\end{align*}
The other inequality is obvious. Let us now show the inequality
\begin{equation}\label{eqnormal}
\d(G_1)\vee G_2 \vee G_3 < G_1 \vee G_2 \vee G_3.
\end{equation}
According to Lemma \ref{lemmasupernormal} it is enough to prove the inequality
\begin{equation}\label{eqsupernormal}
(G_1\vee G_3) \wedge (G_2 \vee G_3) \leq \d(G_1) \vee G_3. 
\end{equation}
An atom $H$ below $G_1 \vee G_3$ is either below $G_1$ or $G_3$ by nested-ness, and similarly for $G_2 \vee G_3$. As a consequence, an atom below $G_1 \vee G_3$ and below $G_2\vee G_3$ is either below $G_3$ or is below $G_1 \wedge G_2$, which is below $\d(G_1)$ by Lemma \ref{lemmasupernormal}. \\

One must also prove that $G_2 \vee G_3$ is not covered by $G_1\vee G_2 \vee G_3$. By inequality \eqref{eqnormal} if $G_2\vee G_3$ is covered by $G_1 \vee G_2 \vee G_3$ then we have $\d(G_1)\vee G_3 \leq G_2 \vee G_3$. In this case by nested-ness $\d(G_1)$ is either below $G_2$ or $G_3$. In the first case we immediately obtain that $G_2$ is covered by $G_1\vee G_2$ which is a contradiction. In the second case we get $G_1 \wedge G_3 \neq \zero $ which contradicts the fact that $G_1\vee G_3$ does not belong to $\G$. \\

\textbf{Case 6.} $G_1 \vee G_2 \in \G$, $G_2 \vee G_3 \in \G$. \\ \\
In this case we necessarily have $G_1 \vee G_2 \vee G_3 \in \G$. Let us prove that $G_1 \vee G_3$ is an initial segment of $G_1 \vee G_2 \vee G_3$ in $[G_3, G_1 \vee G_2 \vee G_3]$. It is enough to prove that $G_1$ is an initial segment of $G_1 \vee G_2 \vee G_3$. By restriction we can assume $G_1 \vee G_2 \vee G_3 = \un$. Since $G_1$ is an initial segment of $G_1 \vee G_2$ there exists an integer $k_1$ such that we have $\d^{k_1}(\un)\wedge (G_1 \vee G_2) = G_1$. Since $G_2$ is an initial segment of $G_2 \vee G_3$ there exists an integer $k_2$ such that we have $\d^{k_2}(\un)\wedge (G_2 \vee G_3) = G_2$. The integer $k_1 $ is greater or equal than $k_2$ because the opposite inequality would imply $G_2 \leq G_1$. Let us prove the equality
\begin{equation*}
\d^{k_1}(\un) = G_1.
\end{equation*}
We have 
\begin{align*}
\d^{k_1}(\un) &= \d^{k_1}(\un)\wedge(G_1 \vee G_2 \vee G_3) \\
&= G_1 \vee (\d^{k_1}(\un) \wedge(G_2 \vee G_3)) \, (\textrm{by modularity of } \d^{k_1}(\un)) \\
&\leq G_1 \vee (\d^{k_1}(\un) \wedge G_2) \, ( \textrm{by } k_1 \geq k_2)\\
&\leq G_1 \vee G_1 \, (\textrm{by definition of } k_1) \\
& = G_1. 
\end{align*}
The opposite inequality is obvious. Let us now prove inequality \ref{eqsupernormal}. By supersolvability the lattice generated by $G_1= \d^{k_1}(\un) , G_3$ and $G_2 \vee G_3$ is distributive. This implies 
\begin{align*}
(G_1 \vee G_3)\wedge (G_2 \vee G_3) &= (G_1 \wedge (G_2 \vee G_3)) \vee (G_3 \wedge (G_2 \vee G_3)) \\
&= (G_1 \wedge (G_2 \vee G_3) ) \vee G_3 \\
&\leq (G_1 \wedge G_2) \vee G_3 \, (\textrm{by } k_1 \geq k_2) \\
&\leq \d(G_1) \vee G_3. 
\end{align*}
As in the other cases one can check that $G_2 \vee G_3$ is not covered by $G_1 \vee G_2 \vee G_3$. \\


\textbf{Case 7.} $G_1 \vee G_2 \in \G$, $G_1 \vee G_3 \in \G$, $G_2 \vee G_3 \notin \G$.\\ \\
In this case we necessarily have $G_1\vee G_2 \vee G_3 \in \G$. As always by restriction we can assume $G_1 \vee G_2 \vee G_3 = \un$. By assumption there exists an integer $k_2$ such that we have $\d^{k_2}(\un) \wedge (G_1 \vee G_2) = G_1$ and an integer $k_3$ such that we have $\d^{k_3}(\un) \wedge (G_1 \vee G_3) = G_1$. Let us denote $k \coloneqq \max(k_2, k_3)$. We will prove the equality
\begin{equation*}
\d^{k}(\un) = G_1. 
\end{equation*}
We have 
\begin{align*}
\d^{k}(\un) &= \d^{k}(\un)\wedge (G_1 \vee G_ 2\vee G_3) \\ 
&= G_1 \vee (\d^k(\un)\wedge(G_2 \vee G_3)) \, (\textrm{by modularity of } \d^k(\un))\\
&= G_1 \vee (\d^k(\un)\wedge G_2) \vee (\d^{k}(\un) \wedge G_3) \, (\textrm{by nested-ness of } \{G_2, G_3\})\\
&\leq G_1 \vee G_1 \, (\textrm{by definition of } k)\\
&= G_1.
\end{align*}
The opposite inequality is obvious. This implies that $G_1 \vee G_3$ is an initial segment of $G_1 \vee G_2 \vee G_3$ in $[G_3, \un]$. Let us now prove inequality \ref{eqsupernormal}. By supersolvability the lattice generated by $G_1= \d^{k_1}(\un) , G_3$ and $G_2 \vee G_3$ is distributive. This implies 
\begin{align*}
(G_1 \vee G_3)\wedge (G_2 \vee G_3) &= (G_1 \wedge (G_2 \vee G_3)) \vee (G_3 \wedge (G_2 \vee G_3)) \\
&= (G_1 \wedge (G_2 \vee G_3) ) \vee G_3 \\
&= (G_1 \wedge G_2)\vee (G_1 \wedge G_3) \vee G_3\, (\textrm{by nested-ness of } \{G_2, G_3\})\\
&\leq \d(G_1) \vee G_3 \,(\textrm{by Lemma \ref{lemmasupernormal}}).
\end{align*}
As in the other cases one can check that $G_2 \vee G_3$ is not covered by $G_1 \vee G_2 \vee G_3$. \\

\end{proof}

\subsection{Proof of the main theorem}
In this subsection we give the proof of our main theorem and some immediate corollaries. 
\begin{theo}
Let $(\L, \G)$ be a supersolvable built lattice. The algebra $\FYa(\L, \G)$ admits a quadratic Gröbner basis.
\end{theo}
\begin{proof}
By dimension it is enough to prove that the map $\Phi$ is a left inverse of $\Psi$, which we will do by induction. The base cases are obvious. \\

Let $\alpha = \prod_{G' \in I} h_{G'}$ be a normal algebraic monomial with maximal element $G$ with respect to $\vartriangleleft$. If $G = \un$ then every element $G' \in I $ is of the form $\d^{k}(\un)$ for some $k$. In this case we can explicitly compute 
\begin{equation*}
\Psi (\prod_{i \leq n} h_{\d^{k_i}(\un)}) = \{\d^{k}(\un) \,  | \, k+1 \notin \{k_i, i \leq n \}\},
\end{equation*}
and 
\begin{equation*}
\Phi (\{\d^{k_i}(\un), \, i \leq n \}) = \prod_{ k-1 \notin \{k_i, i \leq n \}} h_{\d^{k}(\un)},
\end{equation*}
which proves reciprocity. 
If $G \neq \un$ by induction one can prove that the maximal element (for $\vartriangleleft$) of $\Psi(\alpha)$ is $G$. We then have by definition
\begin{align*}
\Phi \circ \Psi(\alpha) &= \Phi(\{G\}\circ(\Psi(G \vee \alpha_{\nleq G}), \Psi(\alpha_{\leq G}))) \\
&= \Supp_{G}(\Phi (\Psi (G\vee \alpha)))\Phi(\Psi(\alpha_{\leq G})). \\
\end{align*}
However, by normality of $\alpha$ and maximality of $G$ we get $\Supp_{G}(G \vee \alpha) = \alpha_{\nleq G}$ which concludes the proof by induction. \end{proof}
We have the immediate corollary.
\begin{coro}
Let $(\L, \G)$ be a supersolvable built lattice. The algebra $\FYa(\L, \G)$ is Koszul. 
\end{coro}
When restricting our attention to the maximal building set we get the following.
\begin{coro}
Let $\L$ be a supersolvable lattice. The combinatorial Chow ring $\FYa(\L, \G_{\max})$ admits a quadratic Gröbner basis. 
\end{coro}
In the next subsection we shall see that this is also true for the minimal building set. 
\subsection{Minimal building sets of supersolvable lattices}
A lattice is said to be irreducible if it is not a product of proper subposets. The minimal building set $\G_{\min}$ of a lattice $\L$ is the set of elements $G$ of $\L$ such that $[\zero, G]$ is irreducible. We have the key proposition. 
\begin{prop}\label{propsupergmin}
Let $\L$ be an irreducible supersolvable lattice. The built lattice $(\L, \G_{\min})$ is supersolvable.
\end{prop}
\begin{proof}
It is enough to prove that if $\L$ is a supersolvable irreducible lattice then $\d(\un)$ is irreducible. In fact here $\d(\un)$ can be any modular coatom (it does not need to be part of a maximal chain of modular elements). We will prove the contraposition of this statement. Assume that $\d(\un)$ decomposes as a product $[\zero, G_1]\times ... \times [\zero, G_n]$. Denote by $R$ the set of atoms which are not below $\d(\un)$. We have the following lemma. 
\begin{lemma}
There exists an integer $i\leq n $ such that we have the inclusion $(\bigvee R) \setminus R \subset G_i$.
\end{lemma}
\begin{proof}
Recall from Subsection \ref{subsecmat} that we have 
\begin{equation*}
\bigvee R = \sigma(R) = R \cup \{i \in \At(\L) \, | \, \exists C \textrm{ circuit s.t. } i \in C \subset R \cup \{i\}\}.
\end{equation*}
Let us prove that there exists some integer $i$ such that for any atom $H \notin R$ and any circuit $C \subset R \cup \{H\}$ containing $H$ we have $H \leq G_i$, by induction on the cardinal of the circuits. The base case is when the circuits have length $3$. Let $H_1$, $H_2$ be two different atoms in $R$. Since $\d(\un)$ is a coatom we have $H_2 \leq \un = \d(\un) \vee H_1$ so there exists a circuit $C$ containing $H_1$ and $H_2$ and such that we have $C\setminus \{H_1, H_2\} \subset \d(\un)$. Since $\d(\un)$ is modular by Lemma \ref{lemmamodularcircuit} there exists a circuit $C'$ containing some element $H'$ in $\d(\un)$ and such that $C'\setminus \{H'\}$ is equal to $\{H_1, H_2\}$. Such a circuit is in fact unique, the atom $H'$ being necessarily equal to $\d(\un)\wedge(H_1 \vee H_2)$. Consider now three atoms $H_1, H_2, H_3$ in $R$. The element $\d(\un) \wedge(H_1 \vee H_2 \vee H_3)$ has rank at most $2$ and it contains the three atoms $\d(\un)\wedge(H_i \vee H_j)$ for $i \neq j \leq 3$. If two of those atoms are equal, say $$ \d(\un) \wedge(H_1 \vee H_2)  = \d(\un) \wedge (H_2 \vee H_3) \eqqcolon H,$$ then by sub-modularity of $\L$ we have $$ \rho(H_1 \vee H_3 \vee H) \leq \rho(H_1 \vee H) + \rho(H_3\vee H) - \rho((H_1\vee H) \wedge (H_3\vee H)) \leq 2 + 2 -2 = 2 $$ (since we have $H \vee H_2  \leq (H_1 \vee H) \wedge (H_3 \vee H)$). This implies $H \leq H_1 \vee H_3$ and therefore $\d(\un) \wedge (H_1 \vee H_3)$ is also equal to $H$. If the three atoms are different then they must form a circuit, and thus they must all belong to some same $G_i$. This concludes the initialization. Let us now assume that all the atoms $\d(\un)\wedge (H_1 \vee H_2)$ are below $G_1$ for instance. \\

Let $C$ be a circuit of arbitrary length, containing a unique element $H$ not in $R$. Let $H_1, H_2$ be two atoms in $C$ different from $H$. By the initialization part there exists a circuit $\{H_1, H_2, H'\}$ with $H'$ in $G_1$. If $H'$ is equal to $H$ then we have $H \in G_1$. If not, by Axiom~\eqref{axcircuit} one can construct a circuit $C'$, containing $H$ and not containing $H_1$, such that $C'$ is included in $C\cup \{H'\}$. If $C'$ does not contain $H'$, then we are done by induction. If $C'$ contains $H'$ then since $\d(\un)$ is modular by Lemma \ref{lemmamodularcircuit} there exists a circuit $C''$ containing some element $H''$ in $\d(\un)$ and such that $C''\setminus\{H''\}$ is contained in $C'\cap R$. By induction the atom $H''$ belongs to $G_1$. If $H'' = H$ we are done. Otherwise by Axiom \eqref{axcircuit} there exists a circuit $C'''$ containing $H$, contained in $C' \cup \{H'' \}$ and not containing some element in $C' \cap R$. Reiterating this process we get a circuit containing $H$ and some elements in $G_1$ which proves that $H$ belongs to $G_1$. 
\end{proof}
From this we deduce the second lemma.
\begin{lemma}
The set of atoms $G_1 \cup R$ is closed. 
\end{lemma}
\begin{proof}
Let $H$ be some element contained in a circuit $C$ contained in $G_1 \cup R \cup \{H\}$. If $H$ is in $R$ we are done. If $H$ is the unique element of $C$ under $\d(\un)$ then by the previous lemma we are also done. Otherwise using Lemma \ref{lemmamodularcircuit} together with Axiom \eqref{axcircuit} (as we did in the previous lemma) gives us a circuit contained in $\d(\un)$, containing $H$ with every other element in $G_1$. This proves that $H$ belongs to $G_1$.
\end{proof}
Finally, we get the concluding lemma. 
\begin{lemma}
$\L$ is isomorphic to $[\zero, G_1\cup R]\times [\zero, G_2] \times ... \times [\zero, G_n]$.
\end{lemma}
\begin{proof}
We will prove that every circuit is either contained in $G_1 \cup R$ or in some $G_i$ with $i \geq 2$. Let $C$ be a circuit in $\L$. If $C$ is contained in $\d(\un)$ then the result comes from the isomorphism $[\zero, \d(\un)] \simeq [\zero, G_1] \times ... \times [\zero, G_n]$. If $C \nsubseteq \d(\un)$ and $C\cap \d(\un)$ is a singleton then by the previous lemma we have $C \subset G_1 \cup R$. If $C\nsubseteq \d(\un)$ and $C \cap \d(\un)$ is not a singleton, pick $H$ any atom in $C\cap \d(\un)$. By iterating Lemma \ref{lemmamodularcircuit} as in the previous proof we obtain a circuit $C'$ containing $H$, contained in $\d(\un)$ and containing some elements in $G_1$. The isomorphism 
\begin{equation*}
[\zero, \d(\un)] \simeq [\zero, G_1]\times [\zero, G_n]
\end{equation*}
implies that this circuit should be contained in $G_1$ which proves the result. 
 \end{proof}
\end{proof}
The above proposition and Theorem \ref{theomain} imply the following theorem.
\begin{theo}
Let $\L$ be a supersolvable lattice. The algebra $\FYa(\L, \G_{\min})$ has a quadratic Gröbner basis and is therefore Koszul. 
\end{theo}
\begin{proof}
Any supersolvable lattice $\L$ decomposes as a product of irreducible supersolvable lattices
\begin{equation*}
\L \simeq \L_1 \times ... \times \L_n.
\end{equation*}
We then have 
\begin{equation*}
\FYa(\L, \G_{\min}) \simeq \FYa(\L_1, \G_{\min})\otimes ... \otimes \FYa(\L_n, \G_{\min})
\end{equation*}
and we can conclude by Proposition \ref{propsupergmin} and Theorem \ref{theomain}.
\end{proof}
\section{Application to the extended modular operad}\label{secgraph}
\subsection{Chordal graphs}
In \cite{Stanley_1972} Stanley proved that if $G$ is a chordal graph (meaning every cycle in $G$ has a chord), the geometric lattice associated to $G$ is supersolvable. This result is based on the following lemma by Dirac \cite{Dirac_1961}.
\begin{lemma}[Dirac, \cite{Dirac_1961}]
Every chordal graph admits a vertex $v$ such that the graph induced by the neighboors of $v$ is a complete graph. 
\end{lemma}
Such vertices are called ``simplicial''. If we remove a simplicial vertex from a chordal graph, the graph we obtain is chordal and this graph is a coatom in the original graph. This means we can reiterate the process and get a maximal chain in the lattice associated with a chordal graph. One can then check that this maximal chain contains only modular elements. We have a ``built'' variant of this result. Let us remind the reader that in Example \ref{exbs} we have defined a built lattice $(\L_G, \G_G)$ for every simple graph $G$, with $\L_G$ the usual graphical matroid associated to $G$ and $\G_G$ the building set of connected subgraphs of $G$. 

\begin{lemma}\label{lemmachords}
Let $G$ be a connected chordal graph. The built lattice $(\L_G, \G_G)$ associated to $G$ is supersolvable.
\end{lemma}
\begin{proof}
Let us choose a maximal chain of modular elements as in the last paragraph. By construction those elements are connected subgraphs of $G$. Let $G'$ be a closed connected subgraph of $G$. For any integer $k$ less than the rank of $G$, the element $\d^{k}(\un) \wedge G'$ can be obtained from $G'$ by successively removing simplicial vertices of $G'$ and therefore it is connected. 
\end{proof}
\begin{coro}\label{corograph}
For all chordal graph $G$, the Feichtner--Yuzvinsky algebra $\FYa(\L_G, \G_G)$ admits a quadratic Gröbner basis. 
\end{coro}

\subsection{The components of the extended modular operad}
In \cite{LM_2000}, Losev and Manin introduced new moduli stacks $\overline{L}_{g,S}$ for stable curves of genus $g$ with painted marked points indexed by $S$ of two types (say ``black'' and ``white'') where the points of type black are allowed to coincide and the points of type white are not. Those stacks are the components of the so-called ``extended modular operad'' (see \cite{LM_2004}). We will deduce from Corollary \ref{corograph} the following result. 
\begin{theo}
The cohomology algebras of the components of the extended modular operads in genus 0 are Koszul. 
\end{theo}
\begin{proof}
It is part of the folklore that if $S$ is a (colored) set with $m$ white points and $n$ black points and $*$ is some chosen white point, then the moduli space $\overline{L}_{0, S}$ is isomorphic to the wonderful compactification of the graphical arrangement $$\{z_i = z_j \, | \,i \neq * \textrm{ white }, \, j \neq * \textrm{ white or black} \},$$ with respect to the building set of connected components (see \ref{exbs}). The corresponding graph denoted $G_{m-1,n}$ has $m+n-1$ vertices, with the first $m-1$ vertices connected to every other vertices and the last $n$ vertices connected only to the first $m-1$ vertices. We notice that $G_{m,n}$ is a chordal graph for every $m$ and $n$ and therefore we can conclude by Corollary~\ref{corograph}. \\

Let us summarize here the main line of arguments giving the stated isomorphism. In the sequel \cite{Manin_2004}, Manin remarked that the moduli stacks $\overline{L}_{g, S}$  are part of the formalism of Hassett spaces introduced by Hassett in \cite{Hassett_2003}. In the latter article, the author introduces the moduli problem of curves with weighted points, where one fixes a ``weight data'' consisting of a vector  $(g, \weight) = (g, w_1, ... , w_n) \in \N \times (]0 , 1] \cap \Q)^n$ and one then seeks to parametrize the nodal curves of genus $g$ with $n$ marked points $(s_i)_{i \leq n}$ which are allowed to coincide ``up to their weights'', meaning that if the points $s_i, i \in I$ coincide then we require
\begin{equation*}  
\sum_{i \in I}w_i \leq 1,
\end{equation*}
and satisfying a (weighted) stability condition (see \cite{Hassett_2003}). If the first $p$ weights are $1$ and the last $n-p$ weights are small enough (precisely $\sum_{i > p} w_i < 1$) then the above condition means exactly that the first $p$ points cannot coincide with any other point and the last $n-p$ points can coincide only between them, and we recover the painted moduli problem of Losev and Manin. Hassett proved that there exists a Deligne--Mumford stack $\overline{\Mod}_{g, \weight}$ representing the above weighted moduli problem. \\

In genus $0$ the stability condition can be simply stated: for any irreducible component $T$ of the nodal curve, we require 
\begin{equation*}
\sum_{i \,\textrm{s.t.}s_i \in T} w_i + \#\textrm{nodes of } T > 2.
\end{equation*}
In addition, in genus $0$ the moduli stack $\overline{\Mod}_{0, \w}$ is a smooth projective scheme (called a Hassett space). If the weights are either $1$ or very small then we call those Hassett spaces ``heavy/light''. Indexes with weight $1$ are called heavy and the other indexes are called light. Work of Cavalieri-Hampe-Markwig-Ranganathan \cite{CHMR_2016} shows that a ``heavy/light'' Hassett space is a tropical compactification of the projective complement $\Mod_{0, \w}$ of the same graphical arrangement $\{z_i = z_j \, | \,i \neq * \textrm{ heavy}, \, j \neq * \textrm{ heavy or light} \}$ (with $*$ some chosen heavy index). To put it in a nutshell this means that there exists an embedding of $\Mod_{0, \w}$ in a torus $\Tor^n$ together with a fan $\Sigma$ in $\R^n$ having support the tropicalization of $\Mod_{0, \w}$ and such that the closure of $\Mod_{0, \w}$ in the toric variety $X(\Sigma)$ is the Hassett space $\overline{\Mod}_{0, \w}$ (we refer to \cite{MS_2015} for an introduction to tropical geometry). The fan $\Sigma$ introduced in \cite{CHMR_2016} is none other than the Bergman fan associated to the buit lattice $(\L_{G_{m,n}}, \G_{G_{m,n}})$ (see \cite{FY_2004} for the definition of the Bergman fan of a built lattice). \\

In \cite{Tevelev_2007}, Tevelev has shown that the tropical compactification of a projective hyperplane arrangement complement along the bergman fan of some building set $\G$ of the correponding lattice can in fact be identified with the wonderful compactification of De Concini and Procesi along the same building set $\G$, which is the stated isomorphism. \end{proof}

\section{Further considerations}\label{secfurther}
\subsection{Towards a classification of Koszul Feichtner--Yuzvinsky algebras}\label{subsecclassify}
We would like to emphasize the fact that we know plenty of Feichtner--Yuzvinsky algebras $\FYa(\L, \G)$  which admit quadratic Gröbner bases and such that the built lattice $(\L, \G)$ is not supersolvable, especially in low rank. For instance, if $C_4$ and $C_5$ are respectively the $4$ and $5$-cycles then the lattices $\L_{C_4}$ and $\L_{C_5}$ are not supersolvable but the built lattices $(\L_{C_4}, \G_{C_4})$ and $(\L_{C_5}, \G_{C_5})$ are so small that they still satisfy the key Lemma \ref{lemmamain} and therefore their Feichtner--Yuzvinsky algebras will be Koszul. However, for the wonderful presentation and for the order used in this article, based on a few examples it feels to the author that the supersolvability condition should be close to necessary, in high enough rank. For instance if $C_n$ is the $n$-cycle then one can easily check that the relations of weight $2$ do not form a Gröbner basis of the algebra $\FYa(\L_{C_n}, \G_{C_n})$ with respect to the order considered in this article for $n \geq 6$. We still do not know if $\FYa(\L_{C_6}, \G_{C_6})$ is Koszul or not. In order to produce a quadratic Gröbner basis of this algebra one would either need to consider a different order, or even a different presentation (which should also be different from the classical presentation since one can show that no order on monomials induces a quadratic Gröbner basis for the classical presentation; the argument is completely analogous to that of Dotsenko \cite{Dotsenko_2022} for the case of the building set of connected subgraphs of the complete graphs).\\

Let us also highlight the fact that even the question of quadraticity of Feichtner--Yuzvinsky algebras is not completely clear. We know that the building sets having a flag nested set complex give quadratic Feichtner--Yuzvinsky algebras but this condition is not necessary, as shown by the following example. Consider $C_4$ the $4$-cycle with edges numbered from $1$ to $4$. The set of flats 
$$ \G = \{\un, \{1,2\}, \{1\}, \{2\}, \{3\}, \{4\}\} $$ is a building set of $\L_{\Circ_4}$ which has a non-flag nested set complex since $\{2, 3, 4\}$ is not nested and does not contain any proper subset which is not nested. However, the Feichtner--Yuzvinsky algebra of this built lattice is the algebra generated by $h_{\un}$ and $h_{\{1,2\}}$ with relations 
\begin{center}
$h_{\un}^3 = 0$, \\
$h_{\un}h_{\{1, 2\}} = h_{\un}^2$, \\
$h_{\{1,2\}}^2 = 0$
\end{center}
which is quadratic since the first relation is a consequence of the last two which are quadratic. This ``pathology'' has to do with the fact that the minimal building set of $\L_{C_4}$ (which is just the atoms together with the maximal element) does not have a flag nested set complex.
\begin{prop}
Let $\L$ be a lattice such that $(\L, \G_{\min})$ has a flag nested set complex. If $\G$ is a building set of $\L$ such that $\FYa(\L, \G)$ is quadratic, then the nested set complex of $(\L, \G)$ is flag.
\end{prop}
\begin{proof}
Assume that we have non-comparable elements $G_1, ... , G_n$, with $n\geq 3$, such that we have $G \coloneqq \bigvee_i G_i \in \G$ and $G_i \vee G_j \notin \G$ for all $i\neq j$. If $G$ is irreducible then decomposing the elements $G_1, ..., G_n$ in their irreducible factors immediately yields a contradiction to the flag-ness of the nested set complex of $(\L, \G_{\min})$. If $G$ is not irreducible, to lighten the notation let us assume $G= \un$ (just restrict to the interval $[\zero, G]$). We have some decomposition
\begin{equation}\label{eqisodec}
\L \simeq [\zero, F_1]\times ... \times [\zero, F_p]
\end{equation}
with irreducible $[\zero, F_i]'s$ and $p \geq 2$. Let $j$ be some index less than $p$. By isomorphism \eqref{eqisodec} we have $F_j = \bigvee_i (F_j \wedge G_i)$. If we decompose the elements $F_j \wedge G_i$ as the join of their factors in $\G$ we can see that there is at most two indexes $i$ such that we have $F_j \wedge G_i \neq \zero$ (otherwise we get a new family of non comparable elements contradicting the flag-ness of $\Ne(\L, \G)$, but this time with join $F_j$ which is irreducible). In addition, there cannot be two such indexes, because if say $F_j = (G_{i_1} \wedge F_j)\vee (G_{i_2} \wedge F_j )$ with $i_1 \neq i_2$ then $F_j$ is an element of $\G$ below $G_1 \vee G_2$ which is neither below $G_{i_1}$ nor below $G_{i_2}$ which contradicts the fact that we have $G_{i_1} \vee G_{i_2} \notin \G$. In conclusion for each $j$ there is exactly one $i$ such that we have $G_i \wedge F_j \neq \zero$, and this implies that we in fact have $F_i \leq G_i$. By using isomorphism \eqref{eqisodec} one more time we get that each $G_i$ is a join of some $F_j$'s, and this forms a partition of the $F_j$'s. Finally, if $\FYa(\L, \G)$ is a quadratic algebra then the relation $(h_{\un}- h_{G_1})...(h_{\un} - h_{G_n})$ can be written as a sum of relations of weight 2, multiplied by monomials. One of the terms of this sum shall be of the form $h_{\un}^{n-2}(h_{\un}- h_{G'_1})(h_{\un} - h_{G'_2})$ with $G'_1$ and $G'_2$ two elements in $\G$ with join $\un$. By isomorphism \eqref{eqisodec} we have $G'_1 = (G_1 \wedge G'_1) \vee ... \vee (G_n \wedge G'_1)$, and similarly for $G'_2$. If there are more than three indexes $i$ such that we have $G_i \wedge G'_1 \neq \zero$, then decomposing the elements $G_i \wedge G'_1$ in $\G$ yields a new obstruction to the flag-ness of $\Ne(\L, \G)$, and we can conclude by some induction. If there are two indexes $i_1\neq i_2$ such that we have $G'_1 \wedge G_{i_1} \neq \zero$ and $G'_1 \wedge G_{i_2} \neq \zero$ then $G'_1$ contradicts the fact that we have $G_{i_1} \vee G_{i_2} \notin \G$. Finally, we get $G'_1 = G_{i _1}$ and $G'_2 = G_{i_2}$ for some $i_1, i_2$, which contradicts $G_{i_1}\vee G_{i_2} \notin \G$. 
\end{proof}
As we know from Proposition \ref{propsupergmin}, if a lattice $\L$ is supersolvable the nested set complex associated to the minimal building set will be flag. 
\subsection{Conceptualizing the proofs of Koszulness}
It would be very beneficial if one could explain in a more conceptual way the strategy for proving the Koszul property introduced by Dotsenko and extended in this paper. In this direction, it could be of interest to check if an analogous strategy could work to reprove the following classical theorem of Yuzvinsky.
\begin{theo}[Yuzvinsky, \cite{Yuzvinsky_2001}]
If $\L$ is a supersolvable geometric lattice then the algebra $OS(\L)$ admits a quadratic Gröbner basis. 
\end{theo}
The corresponding (co)operad would be the cooperad of Orlik--Solomon algebras introduced in \cite{Coron_2022}. Having this other example may lead to a better understanding of the phenomena at play and perhaps give new applications. \\

It would also be interesting to find an operadic characterization of supersolvable lattices, which would explain why they behave so well with respect to the operadic structure. 
\bibliography{QuadGrobBasisBiblio}
\bibliographystyle{plain}
\end{document}